\documentclass[12pt]{article}

\usepackage{hyperref}
\usepackage{amsmath}
\usepackage{amssymb}
\usepackage{amsfonts}
\usepackage{amsopn}
\usepackage{amsthm}
\usepackage[all]{xy}

\swapnumbers
\theoremstyle{plain}
\newtheorem{thm}{Theorem}[section]
\newtheorem{lem}[thm]{Lemma}
\newtheorem{cor}[thm]{Corollary}

\theoremstyle{definition}
\newtheorem{ntt}[thm]{}
\newtheorem{ex}[thm]{Example}
\newtheorem{rem}[thm]{Remark}
\newtheorem{dfn}[thm]{Definition}

\newcommand{\zp}{\mathbb{Z}/p}      

\newcommand{\LL}{\mathbb{L}}   

\newcommand{\DD}{\mathcal{D}}      

\newcommand{\Om}{\Omega}           

\newcommand{\bk}{{\overline{k}}}  
\newcommand{\SX}{X_\bk}           
\newcommand{\SY}{Y_{\bk}}
\newcommand{\SG}{G_{\bk}}         

\newcommand{\res}{\mathrm{res}}     

\DeclareMathOperator{\CH}{\mathrm{CH}}          
\DeclareMathOperator{\Ch}{\mathrm{Ch}}          
\DeclareMathOperator{\CHO}{\overline{\CH}}      

\title{Special correspondences and
Chow traces of Landweber-Novikov operations}

\author{K.~Zainoulline\footnote{Partially supported by SFB 701,
INTAS 05-1000008-8118 and DFG GI~706/1-1.}}

\begin{document}

\maketitle

\begin{abstract}
We prove that the function field of 
a variety which possesses a special correspondence
in the sense of M. Rost preserves rationality of cycles
of small codimensions.
This fact was proven by Vishik in the case of quadrics and
played the crucial role in his construction of fields
with $u$-invariant $2^r+1$.
The main technical tools are the algebraic cobordism of Levine-Morel,
the generalised degree formula and the divisibility of Chow traces
of certain Landweber-Novikov operations.
As a direct application of our methods we prove the similar fact
for all $F_4$-varieties.

\noindent
MSC: 14C15, 55N22
\end{abstract}

\section{Introduction}

In paper \cite{Vi07} A.~Vishik using the techniques of symmetric
operations in algebraic cobordism (see \cite{Vi06})
proved that changing the base field
by the function field of a smooth projective quadric 
doesn't change the property of being rational 
for cycles of small codimension. 
This fact which he calls the Main Tool Lemma
plays the crucial role in his construction of fields
with $u$-invariant $2^r+1$. 

In the present paper we prove
the M.T.L. for a class of varieties introduced by M.~Rost
in the context of the Bloch-Kato conjecture. 
Namely, for varieties 
which possess a {\em special correspondence} (see \cite[Definition~5.1]{Ro06}).
As in Vishik's proof the main technical tools are 
the algebraic cobordism of M.~Levine and F.~Morel, 
the generalised degree formula and 
the divisibility of Chow traces of certain Landweber-Novikov operations.
Therefore, we always assume 
that our base field $k$ has characteristic $0$.

We use the following notation. 
All smooth varieties are assumed to be irreducible.
By $\bk$ we denote the algebraic closure of $k$ and
by $\SX$ the respective base change $X\times_k \bk$ of a variety $X$.
Given a prime $p$ by $\CHO(X)$ 
we denote the Chow ring of $X$ modulo its $p$-torsion part
and by $\Ch(X)=\CHO(X)\otimes\zp$ 
the respective Chow ring with $\zp$-coefficients.
The Chow ring is a graded ring. Its $m$-th graded component
is given by cycles of codimension $m$ and is denoted by $\Ch^m(X)$.
We say that a cycle $y\in\Ch^m(\SX)$ 
is {\it defined over $k$}
if it belongs to the image of the restriction map 
$\res_{\bk/k}\colon\Ch^m(X)\to\Ch^m(\SX)$. 

The following notion will be central in this paper

\begin{dfn} 
Let $X$ be a smooth proper irreducible variety over a field $k$ 
of dimension $n$, $p$ be a prime and $d$ be an integer $0\le d\le n$.
Assume that $X$ has no zero-cycles of degree coprime to $p$.
We say $X$ is a {\it $d$-splitting variety mod $p$} if 
for any smooth quasi-projective variety $Y$ over $k$, for any $m<d$ and
for any cycle $y\in\Ch^m(\SY)$ the following condition holds
\begin{equation}\label{meq}
y\text{ is defined over }k 
\Longleftrightarrow
y_{\bk(X)}\text{ is defined over }k(X). 
\end{equation}
\end{dfn}

\begin{ex}\label{mexam} 
Let $Q$ be an anisotropic projective
quadric over $k$ of dimension $n>2$ and $p=2$.
Then according to A.~Vishik
\begin{enumerate}
\item[(a)] $Q$ is a $[\tfrac{n+1}{2}]$-splitting variety
\cite[Cor.~3.5.(1)]{Vi07}; 
\item[(b)] $Q$ is a $n$-splitting variety
if and only if $Q$ possesses a Rost projector (the proof is unpublished).
\end{enumerate}
\end{ex}

The main result of the paper is the following generalisation 
of \ref{mexam}.(b)

\begin{thm}\label{mthm} Let $X$ be a smooth proper irreducible
variety of dimension $n$ over a field of characteristic $0$.
Assume that $X$ has no zero-cycles of degree coprime to $p$.
If $X$ possesses a special correspondence in the sense of Rost, 
then $X$ is a $\frac{n}{p-1}$-splitting variety 
and the value $\frac{n}{p-1}$ is optimal.
\end{thm}

As an application of the techniques used
in the proof of \ref{mthm}, we provide
a complete list of $d$-splitting 
projective homogeneous varieties of type $F_4$.

\begin{cor}\label{mcor}
Let $X$ be a projective homogeneous variety of type $F_4$ (see Sect~4.IV) 
and $p$ be one of its torsion primes ($2$ or $3$).
Assume that $X$ has no zero-cycles of degree coprime to $p$. Then depending
on $p$ we have
\begin{itemize}
\item[{\footnotesize $p=2$:}] If $X$ is of type $F_4/P_4$, 
then $X$ is a $(\dim X)$-splitting variety.
For all other types $X$ is a $3$-splitting variety and this value is 
optimal.
\item[{\footnotesize $p=3$:}] 
$X$ is always a $4$-splitting variety and this value is optimal.
\end{itemize}
\end{cor}

\begin{ex} An example 
of a non-homogeneous variety which possesses a special correspondence
for $p=3$ was provided recently by N.~Semenov (see \cite{Se08}).
By Theorem~\ref{mthm}
it also provides an example of a $4$-splitting variety.
\end{ex}


\section{Mod-$p$ operations}
In the present section we introduce certain operations 
$$
\phi^{q(t)}_p \colon \Om(X) \to \Ch(X)\text{ parametrised by }
q(t)\in\Ch(X)[[t]]
$$
from the ring of algebraic cobordism $\Om(X)$
to the Chow ring $\Ch(X)$ with $\zp$-coefficients
of a smooth variety $X$, where $p$ is a given prime. 
We also define the Rost number $\eta_p(X)$ and discuss its properties.

\begin{ntt}
The group $\Omega^m(X)$ of cobordism cycles
is generated by classes of proper morphisms
$[Z\to X]$ of pure codimension $m$ with $Z$ smooth.
There are cohomological operations on $\Omega$ parametrised
by partitions 
called Landweber-Novikov operations and denoted by $S_{LN}$.
These operations commute with pull-backs,
satisfy projection and Cartan formulas.  
In our paper we will deal only with operations
given by partitions $(p-1,p-1,\ldots,p-1)$.
Such an operation will be denoted by $S^i_{LN}$, where $i$ is the length
of a partition.

There is a commutative diagram for any integer $m$
$$
\xymatrix{
\Om^m(X) \ar[d]_-{S_{LN}^i} \ar[r]^-{pr} & \CHO^m(X) \ar[r] & 
\Ch^m(X) \ar[d]^-{S^i} \\
\Om^{m+i(p-1)}(X)\ar[r]^-{pr} & \CHO^{m+i(p-1)}(X)\ar[r] & \Ch^{m+i(p-1)}(X)
}
$$
where $S_{LN}^i$ 
is the Landweber-Novikov operation, 
$pr\colon \Om(X)\to \CH(X)$ is the canonical morphism
of oriented theories and
$S^i$ is the $i$-th reduced $p$-power operation.

By properties of reduced power operations $S^i=0$ if $i>m$.
By commutativity of the diagram it means that the composite 
$pr\circ S_{LN}^i$ is divisible by $p$ in $\CHO^{m+i(p-1)}(X)$.
Define (cf. \cite[3.3]{Vi06}) 
\begin{equation}\label{defop}
\phi^{t^{(p-1)a}}_p=\tfrac{1}{p}(pr\circ S_{LN}^{i})\mod p,\text{ where } a=i-m>0.
\end{equation}
If $r$ is not divisible by $(p-1)$, then
we set $\phi^{t^r}_p=0$.
Hence, we have constructed an operation $\phi^{t^r}_p$, $r>0$, 
which maps $\Om^m(X)$ to $\Ch^{r+pm}(X)$.

Finally, given a power series $q(t)\in \Ch(X)[[t]]$ define 
$$
\phi^{q(t)}_p=\sum_{r\ge 0} q_r\phi^{t^r}_p, \text{ where } 
q(t)=\sum_{r\ge 0} q_rt^r.
$$
By the very definition operations $\phi^{q(t)}_p$ 
are additive and respect pull-backs.
\end{ntt}

\begin{dfn}
Let $[U]$ be the class of a smooth projective
variety $U$ of dimension $d$ in the Lazard ring $\LL_d=\Omega^{-d}(k)$. 
Assume that $(p-1)$ divides $d$.
Then the integer $\phi^{t^{dp}}([U])\in \Ch^0(pt)=\zp$
will be called the \emph{Rost number} of $U$ and will be 
denoted by $\eta_p(U)$. 
\end{dfn} 

Using the definition of $S_{LN}$ 
the number $\eta_p(U)$ can be computed as follows.
Let $\xi_1,\xi_2,\ldots,\xi_d$ be the roots of 
the total Chern class of the tangent bundle of $U$.
Define $c(T_U)^{(p)}=\prod_{j=1}^d(1+\xi_j^{p-1})$.
Then 
\begin{equation}\label{etamain}
\eta_p(U)=\tfrac{1}{p}\deg \big(c(T_U)^{(p)}\big)^{-1}.
\end{equation}
Indeed, it coincides with the number $b_d(U)/p$ 
introduced in \cite[Sect.~9]{Ro06}.

\begin{lem}(cf. \cite[Prop.2.3]{Vi07})\label{RDF} 
Let $U$ be a smooth projective variety of positive dimension $d$ 
and $[U]$ be its class in the Lazard ring.
Let $\beta\in \Om^{j}(X)$.
Then
$$
\phi^{t^r}_p([U]\cdot\beta)=\eta_p(U)\cdot pr(S^i_{LN}(\beta)), \text{ where }
r=(p-1)(i-j)+dp>0.
$$
Observe that $\phi^{t^r}_p([U]\cdot\beta)=0$ 
if $d$ is not divisible by $(p-1)$.
\end{lem}

\begin{proof} Let $q\colon X\to Spec\, k$ be the structure map.
Then by Cartan formula
$$
\tfrac{1}{p} (pr\circ S_{LN}^i)([U]\cdot \beta)=\sum_{\alpha+\beta=i}
\tfrac{1}{p} pr (S_{LN}^\alpha(q^*[U])) \cdot 
pr (S_{LN}^\beta(\beta)).
$$
To finish the proof observe that if $\alpha\neq \tfrac{d}{p-1}$, then 
$pr (S_{LN}^\alpha(q^*[U]))=q^* (pr\circ S_{LN}^\alpha([U]))=0$
in $\CHO(X)$ by dimension reasons.
\end{proof}

\begin{cor}(cf. \cite[Lemma~4.4.20]{LM07})
\label{decomcor} Let $U$ and $V$ be smooth projective
varieties of positive dimensions.
Then $\eta_p(U\cdot V)=0$ in $\zp$.
\end{cor}

\begin{proof} Apply Lemma~\ref{RDF} to $X=pt$ and $\beta=[V]$. 
Then $\eta_p(U\cdot V)=\eta_p(U)\cdot 
pr(S^i_{LN}([V]))$, 
where the last factor is divisible by $p$ and, hence, becomes
trivial modulo $p$.
\end{proof}

According to \cite[Remark~4.5.6]{LM07}
the kernel of the canonical morphism $pr\colon \Omega(X)\to\Ch(X)$
is generated by classes of positive dimensions, i.e.
$ker(pr)=\LL_{>0}\cdot\Omega(X)$.
Hence, any $\gamma\in ker(pr)$ can be written as
\begin{equation}\label{kerpr}
\gamma=\sum_{u_Z\in\LL_{>0}} u_Z\cdot [Z\to X].
\end{equation}
Let $\gamma_{pt}\in\LL$ denote the class $u_Z$
corresponding to the point $Z=pt$.

\begin{lem}\label{projpt}
Let $\pi\in\Omega(X\times X)$ be an idempotent and
$\gamma\in \Omega(X)$ be such that
$pr(\pi_\star(\gamma)-\gamma)=0$, where
$\pi_\star$ is the realization.
Then $\eta_p(\pi_\star(\gamma)_{pt})=\eta_p(\gamma_{pt})$.
\end{lem}

\begin{proof}
In presentation \eqref{kerpr}
let $\pi_\star(\gamma)-\gamma=\sum_{u_Z\in\LL_{>0}}u_Z\cdot [Z\to X]$.
Since $\pi$ is an idempotent and $\pi_\star$, $q_*$ are 
$\LL$-module homomorphisms, we obtain
$$
0=q_*\pi_\star(\pi_\star(\gamma)-\gamma)=\sum_{u_Z\in \LL_{>0}} 
u_Z\cdot q_*\pi_\star([Z\to X]).
$$
Apply $\eta_p$ to the both sides of the equality.
By Corollary~\ref{decomcor} all summands with $\dim Z>0$
become trivial (modulo $p$). Hence, $\eta_p(u_{pt})=0\mod p$.
\end{proof}

\begin{cor}\label{idlem} 
Let $\pi$ and $\gamma$ be as above.
Then 
$\eta_p(q_*(\pi_\star(\gamma)-\gamma))=0$, 
where $q\colon X\to Spec\, k$
is the structure map.
\end{cor}

\begin{proof} 
Observe that
$\eta_p(q_*(\pi_\star(\gamma)-\gamma))=
\sum_{u_Z\in\LL_{>0}} \eta_p(u_Z\cdot [Z])$,
where all summands with $\dim Z>0$ are trivial by Corollary~\ref{decomcor}
and the summand with $Z=pt$ is trivial by Lemma~\ref{projpt}.
\end{proof}

The next important lemma is a direct consequence of the result by M.~Rost
\cite[Lemma~9.3]{Ro06}.

\begin{lem}\label{rostm} Let $X$ be a variety which possesses a special correspondence.
Then for any $\gamma\in \Omega_{>0}(X)$ the 
$\deg pr(S^i_{LN}(\gamma))$ is divisible by $p$ in $\CHO(X)$.
\end{lem}

\begin{proof}We have $S^\bullet=S_\bullet\cdot c_\Omega(-T_{\SX})$,
where ($S^\bullet$) $S_\bullet$ are the (co-)homological operations.
Hence, 
$\deg pr(S_{LN}^i(\gamma))=
\deg\big(pr(S_{\bullet}(\gamma))\cdot pr(c_\Omega(-T_{\SX}))\big)$.
Since all Chern classes of the tangent bundle of $\SX$
are defined over $k$ and $X$ possesses a special correspondence,
according to \cite[Lemma~9.3]{Ro06} we obtain that
$$
\deg\big(pr(S_\bullet(\gamma))\cdot pr(c_\Omega(-T_{\SX}))\big)=
\deg pr(S_i(\gamma)) \mod p.
$$
Since $S_\bullet$ respect push-forwards and $\gamma$ has positive dimension,  
$\deg pr(S_i(\gamma))$ is divisible by $p$ as well.
\end{proof}


\section{Construction of a cycle defined over $k$}

In the present section we prove Theorem~\ref{mthm}. The proof
consists of several steps. First, following Vishik's arguments
for a given $y\in\Ch^m(X_\bk)$
we construct a cycle $\bar\omega$ defined over $k$
in the cobordism ring of the product $\Omega(\SX\times\SY)$.
To do this we essentially use the surjectivity of the canonical map 
$pr\colon\Om\to\CH$. 
The motivic decomposition of $X$ provides
an idempotent cycle $\pi$. Applying the realization of $\pi$ to $\bar\omega$
we obtain a cycle $\rho$ defined over $k$ which can be written in the 
form \eqref{mexp}. To finish the proof 
we apply two operations $\phi^{t^r}_p\circ p_{Y*}$
and $p_{Y*}\circ \phi^{t^{r'}}_p$ to the cycle $\rho$.
The direct computations which are based on
the generalised degree formula and \cite[Lemma~9.3]{Ro06} show that
the difference 
$(\phi^{t^r}_p\circ p_{Y*}-p_{Y*}\circ \phi^{t^{r'}}_p)(\rho)$ is defined
over $k$ and provides the cycle $y$.

\paragraph{I.}\label{s3p1}
We start as in the proof of \cite[Thm.~3.1]{Vi07}.
Let $Y$ be a smooth quasi-projective variety over $k$.
Let $y\in\Ch^m(\SY)$ be
such that $y_{\bk(X)}$ is defined over $k(X)$. 
We want to show that $y$ is defined over $k$ for all $m<d$.

Consider the commutative diagram
$$
\xymatrix{
\omega \ar@{|->}[d]& \Om^m(X\times Y) \ar[d]_-{\res}\ar@{>>}[r]^-{pr}& 
\Ch^m(X\times Y)\ar@{>>}[r]^-{pr_1^*}\ar[d]^-{\res} & \Ch^m(Y_{k(X)})\ar[d]^-{\res} 
& u \ar@{|->}[d]\\
\bar\omega &\Om^m(\SX\times \SY) \ar@{>>}[r]^-{pr} & 
\Ch^m(\SX\times\SY)\ar@{>>}[r]^-{\bar pr_1^*} & \Ch^m(Y_{\bk(X)}) & y_{\bk(X)}
}
$$
where the pull-back $pr_1^*$
is surjective by the localisation
sequence and $pr$ is surjective due to \cite[Thm.4.5.1]{LM07}.
By the hypothesis there exists a preimage $u$ of $y_{\bk(X)}$
by means of $\res$. By the surjectivity of $pr$ and $pr_1^*$, there exists
a preimage $\omega$ of $u$.
Set $\bar\omega=\res(\omega)$.

\paragraph{II.}
Let $X$ be a variety which possesses a special correspondence and
has no zero-cycles of degree coprime to $p$.
By the results of M.~Rost it follows that
\begin{itemize}\label{improp}
\item[(a)] $\eta_p(\SX)\neq 0 \mod p$ (see \cite[Thm.~9.9]{Ro06}),
\item[(b)] $n=p^s-1$ (see \cite[Cor.~9.12]{Ro06}), 
\item[(c)] the Chow motive of $X$ contains an indecomposable summand
$M$ which over $\bk$ splits as a direct sum of Tate motives
twisted by the multiples of $d=\tfrac{n}{p-1}$ 
(see \cite[Prop.~7.14]{Ro06})
$$
M_\bk\simeq \bigoplus_{i=0}^{p-1} \zp\{di\}.
$$
\end{itemize}

Let $\pi$ be an idempotent defining the respective $\Omega$-motive $M$. 
Then the realization $\rho=\pi_\star(\bar\omega)$ is defined
over $k$ and 
can be written as (cf.~\cite[p.368]{Vi07})
\begin{equation}\label{mexp}
\rho=x_n\times y_n + \sum_{i=1}^{p-2} 
x_{di}\times y_{di} + x_0\times y_0 \in \Om^m(\SX\times\SY),
\end{equation}
where $x_{j}\in \Om_{j}(\SX)$, 
$y_{j}\in\Om^{m-n+j}(\SY)$,
$x_0=[pt\hookrightarrow \SX]$,
$x_n=\pi_\star(1)$ and $pr(y_n)=y$ (cf.~\cite[Lemma~3.2]{Vi07}).

\paragraph{III.}
Let $p_X^*$ and $p_Y^*$ denote the pull-backs induced by projections 
$X\times Y\to X,Y$.
Since $m<d$, $r=(dp-m)(p-1)>0$.
Consider the cycle $\phi^{t^{r}}_p(p_{Y*}(\rho))$. 
It is defined over $k$ and has codimension $m$. 

\begin{lem}\label{phip} $\phi^{t^r}_p(p_{Y*}(\rho))=
\eta_p(\SX)\cdot y+\phi^{t^r}_p(y_0)$ in $\Ch^m(\SY)$.
\end{lem}

\begin{proof}  
By the projection formula 
$$
\phi^{t^{r}}_p(p_{Y*}(x_{j}\times y_{j}))=
\phi^{t^{r}}_p(p_{Y*}(p_X^*(x_j)\cdot p_Y^*(y_j)))=
\phi^{t^{r}}_p(q_*(x_j)\cdot y_j),
$$
where $x_j\times y_j$
is a summand of \eqref{mexp} and $q\colon X\to Spec\, k$ 
is the structure map.

Assume $0<j<n$. 
By Lemma~\ref{RDF} we obtain
$$
\phi^{t^{r}}_p(p_{Y*}(x_{j}\times y_{j}))=\eta_p(q_*(x_{j}))\cdot 
pr(S^l_{LN}(y_{j})),\text{ where }l=\tfrac{n-j}{p-1}.
$$
Since $codim(y_j)=m-n+j\le m-n+(p-2)d=m-d<0$, we have
$pr(S^l_{LN}(y_i))=S^l(pr(y_j))=0$ in $\Ch^m(\SX)$.
Therefore, only the very right and the left summands of \eqref{mexp} 
remain non-trivial after applying $\phi_p^{t^r}\circ p_{Y*}$.

Now by Lemma~\ref{idlem} the first summand is equal to 
$$
\phi^{t^{r}}_p(p_{Y*}(x_n\times y_n))=\eta_p(q_*(x_n))\cdot pr(y_n)=
\eta_p(q_*(\pi_*(1)))\cdot y=
\eta_p(\SX)\cdot y.
$$
and the last summand
$\phi^{t^{r}}_p(p_{Y*}(x_0\times y_0))=\phi^{t^r}(q_*(x_0)\cdot y_0)=
\phi^{t^r}_p(y_0)$.
\end{proof}

Since $m<d$, $r'=(d-m)(p-1)>0$.
Consider the cycle $p_{Y*}(\phi^{t^{r'}}_p(\rho))$. 
It is defined over $k$ and has codimension $m$.

\begin{lem} \label{pphi}
$p_{Y*}(\phi^{t^{r'}}_p(\rho))=\phi_p^{t^r}(y_0)$ in $\Ch^m(\SY)$.
\end{lem}

\begin{proof} By the very definition
$p_{Y*}\big(\phi^{t^{r'}}_p(x_j\times y_j)\big)=
\tfrac{1}{p}p_{Y*} pr(S_{LN}^{d}(x_j\times y_j))$.
By the projection  and Cartan formulas the latter can be written as
$$ 
\tfrac{1}{p}\deg\big(pr(S_{LN}^{a}(x_j))\big)\cdot pr(S_{LN}^{d-a}(y_j)),
\text{ where }a=\tfrac{j}{p-1}.
$$

Since $m<d$, the cycles 
$y_j$ have negative codimensions for all $j<n$ and, therefore, 
$pr(S_{LN}^{d-a}(y_j))$ is divisible by $p$ for all $j<n$. 
On the other hand, by Lemma~\ref{rostm} the $\deg\big(pr(S_{LN}^{a}(x_j))\big)$
is divisible by $p$ for all $j>0$.
Hence, $p_{Y*}\big(\phi^{t^{r'}}_p(x_j\times y_j)\big)=0$ 
for all $0<j<n$.

Consider the case $j=n$. Recall that $x_n=\pi_\star(1)$. Then we have 
$$
\tfrac{1}{p}\deg\big(pr(S_{LN}^{d}(x_n))\big)=
\sum_{u_Z\in \LL_{>0}}\eta_p(u_Z)\cdot\deg 
pr(S_{LN}^{d-a} ([Z\to X])),
$$
where $a=\tfrac{\dim Z}{p-1}$ and $x_n-1=\sum u_Z[Z\to X]$ in presentation 
\eqref{kerpr}. 
Observe that it is trivial mod $p$, since for all $\dim Z>0$ the degree
of the cycle
$pr(S_{LN}^{d-a} ([Z\to X]))$ is trivial by Lemma~\ref{rostm} and
for $\dim Z=0$ the Rost number $\eta_p(u_{pt})$ is trivial 
 by Lemma~\ref{projpt}.

Hence, only the very last summand, i.e. $x_0\times y_0$,
remains non-trivial after applying $p_{Y*}\circ \phi^{t^{r'}}$
which gives $p_{Y*}(\phi^{t^{r'}}(x_0\times y_0))=\phi^{t^r}(y_0)$.
\end{proof}

\paragraph{IV.}
By Lemmas~\ref{phip} and \ref{pphi} 
the following cycle is defined over $k$
$$
\phi^{t^r}_p(p_{Y*}(\rho))-p_{Y_*}(\phi^{t^{r'}}_p(\rho))=\eta_p(\SX)\cdot y.
$$
Since $\eta_p(\SX)\neq 0\mod p$, the cycle $y$ is defined over $k$, therefore,
$X$ is a $d$-splitting variety.

To see that $d=\tfrac{n}{p-1}$ is an optimal value 
take $Y=X$ and consider the cycle $y\in\Ch^d(\SX)$ 
which generates the Chow group of the Tate motive $\zp\{n-d\}$ 
in the decomposition of $M_\bk$
over $\bk$. Observe that 
$y$ coincides with the cycle $H$
introduced in \cite[Sect.~5]{Ro06}. 
Since $M$ splits over $k(X)$, $y_{\bk(X)}$ is defined over $k(X)$.
By condition \eqref{meq} we obtain that $y$ is defined over $k$, i.e.
$M$ splits over $k$ which contradicts to the indecomposability of $M$.
The theorem is proven.

\section{$F_4$-varieties}

In the present section we apply our methods to describe
all $d$-splitting varieties of type $F_4$.

\paragraph{I.}
Let $X$ be a smooth geometrically cellular 
variety over $k$ of dimension $n$. 
As in the beginning of the previous section
given a cycle $y\in\Ch^m(\SY)$, where $Y$ is smooth quasi-projective,
we construct a cobordism cycle $\bar\omega\in\Om^m(\SX\times\SY)$ 
defined over $k$. 

\paragraph{II.} 
Assume that the motive of $X$ contains a motive $M=(X,\pi)$ 
such that
$$
\pi_\bk=\gamma\times \gamma^\vee + \sum_{i=0}^{p-2} x_{di}\times x_{di}^\vee,
\text{ where }d>1,
$$
the cycle $\gamma\in\Om_{(p-1)d}(\SX)$ is defined over $k$, 
$\gamma^\vee$ denotes its Poincare dual, i.e. 
$pr(\gamma\cdot\gamma^\vee)=pt$, and
$x_j\in \Om_j(\SX)$.
Then the realization
$\rho=\pi_\star(p_X^*(\gamma)\cdot\bar\omega) \in\Om^{m+n-g}(\SX\times\SY)$ 
is defined over $k$
and can be written as (cf.~\eqref{mexp})
$$
\rho=x_g \times y_g + \sum_{i=1}^{p-2} x_{di} \times y_{di}+ 
x_0\times y_0,
$$
where $g=(p-1)d$,
$x_g=\pi_\star(\gamma)$ and $pr(y_g)=y$.

The transposed cycle $\pi^t$ defines an opposite 
direct summand $M^t=(X,\pi^t)$
of the motive of $X$ 
(the one which contains the generic point of $X$ over $\bk$).
The realization $\rho'=\pi^t_\star(\bar\omega)$
is defined over $k$ and can be written as
$$
\rho'=x^{(0)}\times y^{(0)}+\sum_{i=1}^{p-2} 
x^{(di)}\times y^{(di)} + x^{(g)}\times y^{(g)} 
\in\Om^m(\SX\times\SY),
$$
where $x^{(j)}\in\Om^j(\SX)$ and 
$y^{(g)}=y_0$.

\paragraph{III.} 
Since the Chow group
of the cellular variety $X_\bk$ is torsion-free,
we may use Mod-$p$ operations over $\bk$.
We now apply the operations $\phi^{t^r}\circ p_{Y*}$ and
$p_{Y*}\circ(p_X^*(\gamma)\cdot \phi^{t^{r'}})$ to the cycles $\rho$
and $\rho'$ respectively.

Applying the first operation and
repeating the arguments of the proof of Lemma~\ref{phip} we obtain that 
for any $m<d$, $r=(dp-m)(p-1)$,   
the following cycle is defined over $k$
$$
\phi^{t^r}_p(p_{Y*}(\rho))=\eta_p(\gamma)\cdot y+
\phi^{t^r}_p(y_0)\text{ in }\Ch^m(\SY).
$$

To apply the second operation 
for any $m<d$ consider the cycle $\delta=p_{Y*}(p_X^*(\gamma)\cdot\phi^{t^{r'}}_p(\rho'))$,
where $r'=(d-m)(p-1)$.
It is defined over $k$ and has codimension $m$.
Since we don't have the version of \cite[Lemma~9.3]{Ro06}
for an arbitrary variety, to compute $\delta$ we have
to treat each torsion prime case separately.

For $p=2$ we obtain $\delta=\phi^{t^r}_p(y_0)$
by dimension reasons. Indeed, in this case the cycle $\rho'$
consists only of two summands
$$
\rho'=x^{(0)}\times y^{(0)}+\gamma^\vee\times y^{(g)},
$$ 
where the first summand vanishes, 
since $S^{\alpha}_{LN}(x^{(0)})=0$ if
$|\alpha|>0$ and the second summand gives the required cycle 
$\phi^{t^r}_p(y_0)$.

For $p=3$ the cycle $\rho'$ consists of three terms
$$
\rho'=x^{(0)}\times y^{(0)}+x^{(d)}\times y^{(d)}+ \gamma^\vee\times y^{(g)},
$$
where again the first summand vanishes, 
the last gives $\phi^{t^r}(y_0)$
and the middle gives
\begin{equation}\label{mids}
\deg\big(pr(\gamma)\cdot pr(S^{\alpha}_{LN}(x^{(d)}))\big)
\cdot 
\tfrac{1}{p}pr(S^{\beta}_{LN}(y^{(d)})),\text{ where } 
\alpha=\beta=d/2.
\end{equation}

Hence, following the part IV of the previous section
to prove that $X$ is a $d$-splitting variety for $p=2$ or $3$
it is enough to assume that $\eta_p(\gamma)\neq 0\mod p$ and 
that \eqref{mids} vanishes for $p=3$.

\paragraph{IV.}
We use the following notation.
Let $G$ be a simple linear algebraic group over $k$.
We say a projective homogeneous $G$-variety $X$ 
is of type $\DD$, 
if the group $\SG$ has a root system of type
$\DD$. Moreover, if 
$\SX$ is the variety of parabolic subgroups of $\SG$
defined by the subset of simple roots $S$ of $\DD$, then
we say that $X$ is of type $\DD/P_S$.
In this notation $P_\DD$ defines a Borel subgroup
and $P_i$, $i\in\DD$, defines a maximal parabolic subgroup
(our enumeration of roots follows Bourbaki). 

Given an $F_4$-variety $X$ we provide a cycle $\gamma$ satisfying
$\eta_p(\gamma)\neq 0\mod p$ as follows

\subparagraph{$p=2$:} If $X$ is generically split 
over the 2-primary closure of $k$, 
then we may assume that $X$ is of type $F_4/P_1$.
In this case $X$ has dimension $15$ and 
by the main result of \cite{PSZ} the Chow motive
of $X$ with $\zp$-coefficients splits as a direct sum
of twisted copies of a certain motive $M=(X,\pi)$ 
with the generating function
$P(M_\bk,t)=1+t^3$ (see p.33 of \cite{PSZ}).
Since the Chow group $\Ch_r(\SX)$ has rank $1$ for $r=0\ldots 3$,
the idempotent $\pi_\bk$ can be written as
$\pi_\bk=\gamma \times\gamma^\vee + pt \times 1$, where 
$\gamma$ is represented by a 3-dimensional subquadric
$Q_3\hookrightarrow \SX$ which is 
an additive generator of $\Ch_3(\SX)$ defined over $k$.
Since $\eta_2(Q_3)=1$, $X$ is a $d$-splitting variety with $d=3$.

If $X$ is not generically split, i.e. $X$ is of the type $F_4/P_4$,
then $X$ is a splitting variety of the symbol given by the cohomological
invariant $f_5$. By the result of Rost \cite[Rem.~2.3 and \S~8]{Ro06}
$X$ is a variety which possesses a special correspondence.
Hence, by Thm.~\ref{mthm} $X$ is a $15$-splitting variety.

\subparagraph{$p=3$:} 
In this case all $F_4$-varieties are generically split over the $3$-primary
closure of $k$
and we may assume that $X$ is of type $F_4/P_4$.
Similar to the previous case using the motivic
decomposition from \cite{PSZ}
we obtain an idempotent
$$
\pi_\bk=\gamma\times\gamma^\vee+ x_4\times x_4^\vee + pt\times 1,
$$ 
where $\gamma\in\Omega_8(\SX)$ 
is defined over $k$ and $x_4\in\Omega_4(\SX)$ is not.
By the explicit formulae from \cite[5.5]{NSZ} we may identify $pr(\gamma)$
with the $7$-th power of the generator $H$ of the Picard group of $\SX$. 
which is the only cycle in $\Ch_8(\SX)$ defined over $k$.
Since $H$ is very ample, $H^7$ can be 
represented by a smooth projective subvariety $Z$
of $\SX$. Hence, we may identify $\gamma$ with the class
$[Z\hookrightarrow \SX]$.
Then the direct computations using 
the adjunction formula \cite[Example~3.2.12]{Fu} 
show that
$\eta_3(Z)\neq 0\mod 3$.  

To prove the vanishing of the cycle \eqref{mids} it is enough to prove
the vanishing of the cycle
$$
\deg\big(pr(\gamma)\cdot S^2(pr(x^{(4)}))\big)\cdot \phi^{t^{12-2m}}_3(y^{(4)}).
$$
Direct computations show that $S^2(\Ch^4(\SX))$ is trivial, hence,
$\delta=\phi^{t^r}_p(y_0)$ and
$X$ is a $d$-splitting variety with $d=4$. 
The Corollary~\ref{mcor} is proven.

\begin{rem} Let $X$ be a $d$-splitting geometrically cellular
variety. As an immediate consequence of \cite[Cor.~4.11]{KM06} we obtain 
the following bound for a canonical $p$-dimension of $X$
$$
cd_p(X)\ge d.
$$
In the case of a variety of type $F_4/P_4$ it gives
$cd_2(X)=\dim X=15$.
\end{rem}


\paragraph{Acknowledgements} I am very grateful to Alexander Vishik and Fabien
Morel for various discussions on the subject of this paper.

\bibliographystyle{chicago}

{\small \noindent
K.~Zainoulline, Mathematisches Institut der LMU M\"unchen,\\ 
Theresienstr.~39,
D-80333 M\"unchen}

\end{document}